\documentclass{amsart}
\usepackage[frame,web]{xy}
\usepackage{url,amsmath,amssymb}

\newcommand{\abs}[1]{\left\lvert #1\right\rvert}
\let\size=\abs
\newcommand{\eval}[2][\right]{\relax
  \ifx#1\right\relax \left.\fi#2#1\rvert}
\newcommand{\Set}[1]{\ensuremath{\mathcal{#1}}}            
\newcommand{\Dfn}[1]{\emph{#1}}                            
\newtheorem{thm}{Theorem}[section]

\newtheorem{cnj}[thm]{Conjecture}

\DeclareMathOperator{\length}{length}
\DeclareMathOperator{\xfinal}{x-final}
\DeclareMathOperator{\yfinal}{y-final}
\DeclareMathOperator{\res}{res}
\newcommand{\BigO}{O}
\newcommand{\littleo}{o}

\author{Martin Rubey}
\begin{document}

\title{Transcendence of generating functions of walks on the slit plane}
 \thanks{LaBRI, Universit\'e Bordeaux I, Research financed by
   EC's IHRP Programme, within the Research Training Network \lq\lq Algebraic
   Combinatorics in Europe\rq\rq, grant HPRN-CT-2001-00272.
   \url{martin.rubey@labri.fr}\\
   \url{http://www.mat.univie.ac.at/~rubey}}

\begin{abstract}
  Consider a single walker on the slit plane, that is, the square grid $\mathbb
  Z^2$ without its negative $x$-axis, who starts at the origin and takes his
  steps from a given set $\mathfrak S$. Mireille Bousquet-M\'elou conjectured
  that -- excluding pathological cases -- the generating function counting the
  number of possible walks is algebraic if and only if the walker cannot cross
  the negative $x$-axis without touching it. In this paper we prove a special
  case of her conjecture.
\end{abstract}
\maketitle

\section{Introduction}
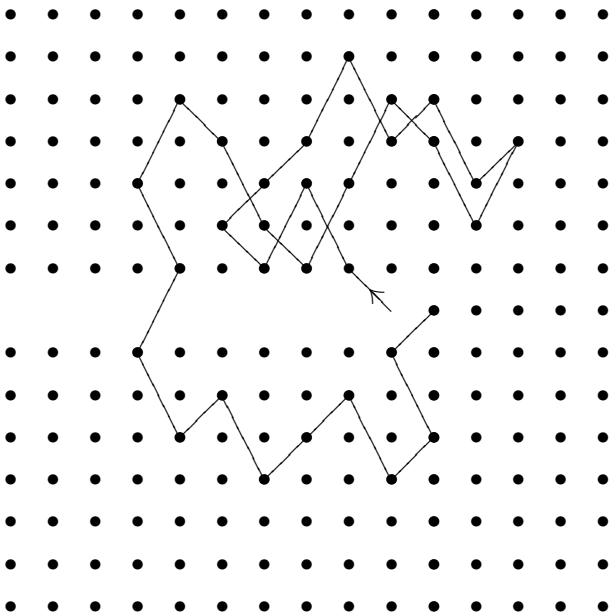
\begin{figure}[t]
  \begin{equation*}
  \begin{xy}<16pt,0pt>:
    {\xylattice{0}{4}{-7}{7}\xylattice{-10}{0}{1}{7}\xylattice{-10}{0}{-7}{-1}}
    @i @={(-1,2),(-1,-2),(-1,1),(1,1),(1,1),(1,2),(1,-2),(1,1),(1,-2),
          (1,1),(-1,-2),(-1,2),(-1,1),(-1,-2),(-1,-2),(-1,1),(-1,2),(-1,1),
          (-1,-2),(1,-2),(-1,-2),(1,-2),(1,1),(1,-2),(1,1),(1,1),(1,-2),(1,1),
          (-1,2),(1,1)},
    (-1,1)="prev", @@{;p+"prev";"prev";**@{-}="prev"*{\bullet}},
    (0,0);(-1,1)**\dir{-}?*\dir3{>}
  \end{xy}
  \end{equation*}
  \caption{\label{fig:walks}A walk on the slit plane}
\end{figure}

Let $\mathfrak S$ -- the set of \Dfn{steps} -- be a finite subset of $\mathbb
Z^2$. A \Dfn{walk on the slit plane} is a sequence $(0,0)=w_0,w_1,\dots,w_n$ of
points in $\mathbb Z^2$, such that the difference of two consecutive points
$w_{i+1}-w_i$ belongs to the set of steps $\mathfrak S$ and none of the points
but the first lie on the half-line $\{(x,0):x\leq 0\}$. An example
for such a walk with set of steps
$$\mathfrak S=\{(-1,-2),(-1,1),(-1,2),(1,-2),(1,1),(1,2)\}$$
is shown in Figure~\ref{fig:walks}.

Recall that a generating function $F(t)=\sum_{n\ge 0} f_n t^n$ is
\Dfn{algebraic}, if there is a nontrivial polynomial $P$ in two variables, such
that $P(F(t),t)=0$. Otherwise, it is \Dfn{transcendental}.

In \cite{BousquetMelou2001} Mireille Bousquet-M\'elou conjectured the
following:
\begin{cnj}\label{cnj:Mireille}
  Consider the generating function for walks in the slit plane with a given set
  of steps $\mathfrak S$, counted according to their length and their
  end-coordinates:
  $$S(x,y;t)=\sum_{\substack{\text{$W$ walk on the slit plane}\\ 
                             \text{starting at the origin}\\ 
                             \text{with steps in $\mathfrak S$}}}
  t^{\length W} x^{\xfinal W} y^{\yfinal W}.$$
  Suppose that the set of steps is not degenerated and thus all four quadrants
  of the plane can be reached by some walk, and that the greatest common
  divisor of the vertical parts of the steps is equal to one.
  
  Then this generating function is algebraic in $t$, if and only if the height
  of any step is at most one.
\end{cnj}

In fact, she proved one part of this conjecture in Section~7 of the above
paper, namely, that walks with steps that have height at most one have an
algebraic generating function. Furthermore, in Section~8 she proved for one
family of step-sets that the corresponding generating functions have to be
transcendental. In the present paper, we prove the following:

\begin{thm}\label{thm:transcendental}
  Let \Set H and \Set V be two finite sets of integers, the greatest common
  divisor of the integers in each set being equal to one. Furthermore, assume
  that both of the sets \Set H and \Set V contain positive and negative
  numbers, and that \Set V contains an element with absolute value at least
  $2$. Finally, assume that the minimum of \Set V is at least $-2$.
  
  Let $\mathfrak S$ be the Cartesian product of the two sets: $\mathfrak S=\Set
  H\times\Set V$, where \Set H is the horizontal and \Set V is the vertical
  part of the steps. Then the following generating functions for walks in the
  slit plane with set of steps $\mathfrak S$ are transcendental in $t$: 
  \begin{itemize}
  \item the generating function $S_{i,0}(t)$ for walks ending at a prescribed
    coordinate $(i,0)$,
  \item the generating function $L(t)$ for loops, i.e., walks that return to
    the origin,
  \item the generating function $S_0(1;t)$ for walks ending anywhere on the
    x-axis, and
  \item $S(1,1;t)$, which is the generating function for walks ending anywhere
    in the slit plane.
  \end{itemize}
\end{thm}

For example, the set of steps of the walk in Figure~\ref{fig:walks} is the
Cartesian product of $\Set H=\{-1,+1\}$ and $\Set V=\{-2,+1,+2\}$.

In fact we consider a slightly more general problem: we allow the steps in
$\Set H$ and $\Set V$ to be weighted with positive real numbers. The weight of
a step in the product set $\mathfrak S=\Set H\times\Set V$ then is the product
of the weights of its corresponding vertical and horizontal parts and the
weight of a walk is the product of the weights of its individual steps.

As in \cite{BousquetMelou2001}, we will use a special case of the following
theorem to determine in which cases the generating function cannot be
algebraic:
\begin{thm}\cite{Flajolet1987}
  Let $F(t)$ be an algebraic function over $\mathbb Q$ that is analytic at the
  origin, then its $n$\textsuperscript{th} Taylor coefficient $f_n$ has an
  asymptotic equivalent of the form
  $$ f_n=\frac{\beta^n n^s}{\Gamma(s+1)}
         \sum_{i=0}^m C_i\omega_i^n+\BigO(\beta^n n^t),$$
  where $s\in\mathbb Q\setminus\{-1,-2,-3,\dots\}$, $t<s$; $\beta$ is a
  positive algebraic number and the $C_i$ and $\omega_i$ are algebraic with
  $\abs{\omega_i}=1$. 
\end{thm}

It follows easily that an algebraic function cannot have an appearance of a
negative integer power of $n$ anywhere in the full asymptotic expansion of its
Taylor coefficients.

\section{An expression for the generating function for walks on the slit plane}

The fundamental theorem for walks on the slit plane is the following:
\begin{thm}(Proposition 9 in \cite{BousquetMelou2001})\label{thm:cycle}
  Let 
  $$B(x;t)=\sum_{\substack{\text{$W$ walk on $\mathbb Z^2$}\\
                           \text{starting at the origin}\\
                           \text{ending on the $x$-axis}\\
                           \text{with steps in $\mathfrak S$}}}
  x^{\xfinal(W)}t^{\length(W)}$$
  be the generating function for bilateral
  walks, that is, walks that end on the $x$-axis but are otherwise
  unconstrained.

  For $i\ge 1$, the generating function $S_{i,0}(t)$ for walks on the slit
  plane ending at $(i,0)$ can be computed by induction on $i$ via the following
  identity: 
  $$
  \sum_{k=1}^i\frac{(-1)^{k-1}}{k}\sum_{\substack{i_1+i_2+\dots+i_k=i\\
      i_1>0,i_2>0,\dots,i_k>0}} S_{i_1,0}(t)S_{i_2,0}(t)\dots S_{i_k,0}(t)
   = [x^i]\log B(x;t).
  $$
\end{thm}
Note that it follows that $\log B(x;t)$ has positive Taylor coefficients.

Now we can take advantage of the special structure of the set step $\mathfrak
S$. Since it decomposes into a horizontal and a vertical part, the generating
function for bilateral walks factorises:
$$ B(x;t)=\bar B(H(x)t)$$
where
$$ \bar B(t)=\sum_{\substack{\text{$W$ walk on $\mathbb Z$}\\
                               \text{from $0$ to $0$}\\
                               \text{with steps in $\Set V$}}}
                           t^{\length(W)}
$$
is the generating function for \Dfn{bridges} and
$$ H(x)=\sum_{h\in\Set H}x^h$$
is the \Dfn{step (Laurent-)polynomial} for $\Set H$.

\section{Transcendence of $[x^i]\log B(x;t)$}
In this section we show that $[x^i]\log B(x;t)$ cannot be algebraic if the
set $\Set V$ contains an element with absolute value strictly greater than one.

To this end, we consider the asymptotic expansion of $[t^n x^i]\log
B(x;t)$. Since $B(x;t)$ factorises, we have
\begin{equation}\label{eq:factorisation}
  \begin{split}
    [t^n x^i] \log B(x;t)&=[t^n x^i] \log \bar B\left(H(x)t\right)\\
    &=[x^i] \left(H(x)\right)^n [t^n]\log \bar B(t).
  \end{split}
\end{equation}

Therefore, we have divided the problem in two: we will show that the asymptotic
expansions of both $[x^i]\left(H(x)\right)^n$ and $[t^n]\log \bar B(t)$ contain
a term $n^{-k/2}$ for some odd $k$.

\subsection{Asymptotics of the horizontal part}
Let $min$ be the minimal integer such that $H(x)x^{min}$ is a polynomial. To
determine the asymptotics of $[x^i] \left(H(x)\right)^n=[x^{n\cdot min}]x^{-i}
\left(H(x)x^{min}\right)^n$, we can use the following theorem:
\begin{thm}[\cite{Drmota1994,FlajoletSedgewick6}]
  Let $g(z)$ be an analytic function of degree $d$ with positive
  coefficients assumed to be aperiodic and such that $g(0)\neq 0$, and let
  $a(z)$ be analytic except possibly at zero, where a pole is allowed. Let
  $\lambda$ be a positive number of some subinterval $[\lambda_a,\lambda_b]$ of
  the open interval $]0,d[$. Then, with $N=\lfloor\lambda n\rfloor$, one has
  uniformly for $\lambda \in [\lambda_a,\lambda_b]$
  $$ [z^N]a(z)(g(z))^n= a(\zeta)\frac{(g(\zeta))^n}{\zeta^{N+1}\sqrt{2\pi n
      R}}(1+\littleo(1)),$$
  where $\zeta$ is the unique positive root of the equation
  $$ \zeta\frac{g^\prime(\zeta)}{g(\zeta)}=\lambda$$
  and
  $$ R=\frac{d^2}{d\zeta^2}[\log g(\zeta)-\lambda\log\zeta].$$
\end{thm}

Let $p=\gcd\{k+min:k\in\Set H\}$ and note that in the case $p>1$, every
$p$\textsuperscript{th} coefficient of $\left(H(x)\right)^n$ will be zero.
Thus, we take $g(x)=H(x^{1/p})x^{min/p}$, and $a(x)=x^{-i/p}$, and use the
theorem to determine $[x^{n\cdot min/p}]a(x)(g(x))^n$. In this situation,
$\lambda=min/p$ is constant, therefore $\zeta$ and $R$ must be constant, too.
Hence, the asymptotic expansion of $[x^i] \left(H(x)\right)^n$ contains an
appearance of $n^{-1/2}$.

\subsection{Asymptotics of the vertical part}
To determine the asymptotic behaviour of $[t^n]\log \bar B(t)$ we will use
analysis of singularities. In a first step, we have to determine the
singularities of the expression. Following the general theory of singularity
analysis, all the contributions from these singularities must be added up.


According to the remark after Theorem~\ref{thm:cycle} and the
factorisation~\eqref{eq:factorisation}, we have that $[t^n]\log \bar B(t)$ is
always positive. Therefore we can apply Pringsheim's theorem:
\begin{thm}\cite{Hille}
  If a function with a finite radius of convergence has Taylor coefficients
  that are nonnegative, then one of its singularities of smallest modulus -- a
  dominant singularity -- is real positive.
\end{thm}

Since the logarithm is singular only at the origin, and $\bar B(t)$ is strictly
positive for positive real numbers, Pringsheim's theorem implies that one
dominant singularity of $\log \bar B(t)$ is in fact a singularity -- call it
$\rho$ -- of $\bar B(t)$. Of course, there can be other dominant singularities
of $\log \bar B(t)$ arising from singularities of $\bar B(t)$. We will discuss
them below. 

Furthermore, it might happen -- although we believe that it does not -- that
$\bar B(t)$ has zeros on the circle around the origin with radius $\rho$,
thus also making $\log\bar B(t)$ singular. However, such singularities can
never contribute a summand of order $n^{-k/2}$ for odd $k$, so we can simply
ignore them. Note that $\bar B(t)$ cannot vanish for $\abs{t}<\rho$, since
$\rho$ is a dominant singularity of $\log \bar B(t)$.

Now we want to compute the contribution of the dominant singularities of $\bar
B(t)$ to the asymptotic expansion of $\log\bar B(t)$. To do so, we need a
better understanding of $\bar B(t)$, which is the generating function for
bridges with step set $\Set V$.

Luckily, this generating function has already been studied. We define the
\Dfn{step (Laurent-)polynomial} for $\Set V$ as
\begin{equation*}
  V(y)=\sum_{v\in\Set V}y^v
\end{equation*}
and the \Dfn{characteristic curve} determined by $\Set V$ by the equation
\begin{equation}
  \label{eq:characteristiccurve}
  1-t V(y)=0\quad\text{or equivalently}\quad y^{min_v}=t(y^{min_v}V(y))=0,
\end{equation}
where $min_v=-\min \Set V$ is the minimal integer to make the equation
polynomial.

We say that the functional equation~\eqref{eq:characteristiccurve} is
\Dfn{reduced}, if the greatest common divisor of the exponents of the monomials
in $V(y)$ is equal to one, which is one of the assumptions in our main
Theorem~\ref{thm:transcendental}.

We say that the functional equation~\eqref{eq:characteristiccurve} has
\Dfn{period} $p$, if the greatest common divisor of the exponents of the
monomials in $y^{min_v}V(y)$ is equal to $p$.

As is well known, the period is also the number of dominant singularities of
$\bar B(t)$, which are all conjugate to the real dominant singularity $\rho$.
In our case however, it can be seen (\cite[Section~3.3]{BanderierFlajolet2002})
that the asymptotic formula for bridges is obtained from the asymptotic
expansion derived from the singularity at $\rho$ by multiplying with $p$. Since
we are only interested in the presence or absence of a term $n^{-k/2}$ for some
odd $k$ in the asymptotic expansion, we can assume from now on that the
functional equation~\eqref{eq:characteristiccurve} is aperiodic, i.e., has
period one.

It can be seen \cite{BanderierFlajolet2002,Hille} that the solutions of this
functional equation organise themselves into \lq\lq small\rq\rq\ and \lq\lq
large\rq\rq\ branches. Here, \lq\lq small\rq\rq\ means that the solution $y(t)$
tends to zero as $t$ tends to zero, whereas \lq\lq large\rq\rq\ means that
$y(t)$ tends to infinity as $t$ approaches zero.

It is only the set of \lq\lq small\rq\rq\ solutions that is interesting for us,
and it can be seen -- using a limit case of Pellet's Theorem, see for example
\cite{Marden} -- that there are $min_v$ of them. A nice expression for the
generating function for bridges is given by the following theorem:
\begin{thm}\cite[Theorem~1 and proof of Theorem~3]{BanderierFlajolet2002}\label{thm:bridges}
  The generating function for bridges is an algebraic function given by

  \begin{equation*}
    \begin{split}
      \bar B(t)&=t\sum_{j=1}^{min_v}\frac{y_j^\prime(t)}{y_j(t)}
                =t\frac{d}{dt}\log(y_1(t)y_2(t)\cdots y_{min_v}(t))\\
               &=\frac{1}{2\pi i}\int_{\abs{y}=\tau}\frac{dy}{y(1-tV(y))},
    \end{split}
  \end{equation*}
  where the expressions involve all the small branches
  $y_1,y_2,\dots,y_{min_v}$ of the characteristic
  curve~\eqref{eq:characteristiccurve}, and $\rho$ is the radius of
  convergence of $\Bar B(t)$. Furthermore, the principal branch
  $y_1(t)$, i.e., the branch with real coefficients, has a square root
  singularity at $\rho$ and the product of all the other small
  branches is analytic for all $t$ with $\size t \leq\rho$. More
  precisely, $\rho$ is given by $\rho=1/V(\tau)$, where $\tau$ is the
  unique positive number with $V^\prime(\tau)=0$.
\end{thm}


Thus, applying the Newton-Puiseux theorem, we can develop $\bar B(t)$ around
the singularity $\rho$, setting $\tilde t=\sqrt{\rho-t}$:
\begin{equation*}
  \bar B(t)=a_{-1}/\tilde t + a_0 + a_1 \tilde t + a_2 \tilde t^2+\cdots
\end{equation*}

Composing this expansion with the Taylor expansion of the logarithm we obtain
\begin{equation*}
  \log\bar B(t)=\log a_{-1}/\tilde t +
  \log(1+\frac{a_0}{a_{-1}}\tilde t+\frac{a_1}{a_{-1}}\tilde t^2+\cdots). 
\end{equation*}

Now we want to find a term $n^{-k/2}$ for some odd $k$ in the
asymptotic expansion of the coefficients of the above series. Clearly,
\begin{equation*}
  [t^n]\log a_{-1}/\tilde t=[t^n]\log a_{-1}/\sqrt{\rho-t}\sim c_0\rho^{-n}n^{-1},
\end{equation*}
for some constant $c_0$, is not what we are looking after. However, we have
\begin{equation*}
  \log(1+\frac{a_0}{a_{-1}}\tilde t+\frac{a_1}{a_{-1}}\tilde t^2+\cdots)
  =\frac{a_0}{a_{-1}}\tilde t+\cdots
\end{equation*}
and
\begin{equation*}
  [t^n]\frac{a_0}{a_{-1}}\tilde t=[t^n]\frac{a_0}{a_{-1}}\sqrt{\rho-t}\sim c_1\rho^{-n}n^{-3/2},
\end{equation*}
for some constant $c_1$.  Provided that $a_0$ does not vanish, this term will
guarantee transcendence of the generating function for walks on the slit plane.
Thus we need the constant term in the singular expansion of $\bar B(t)$.

Since the product of the non-principal branches $y_2(t)y_3(t)\dots y_{min_v}(t)$
is analytic and non-zero at $\rho$, the contribution of
$t\frac{d}{dt}\log(y_2(t)\cdots y_{min_v}(t))$ to the constant term in the
singular expansion of $\bar B(t)$ around $\rho$ is the sum of the residues of
$1/y(1-\rho V(y))$ at the zeros of $1-\rho V(y)$ that are strictly smaller than
$\tau$ in modulus.

To obtain the contribution of $t\frac{d}{dt}\log(y_1(t))$, we proceed
as follows:
\begin{equation*}
  \begin{split}
    &[\tilde t^0]t\frac{d}{dt}\log(y_1(t))\\ 
    &=[\tilde t^0](\rho-\tilde t^2)(-\frac{1}{2\tilde t})\frac{d}{d\tilde t}
    \log(y_1(t))\\
    &=-\frac{\rho}{2}[\tilde t]\frac{d}{d\tilde t}\log(y_1(t))\\
    &=\rho[\tilde t^2]\log(y_1(t)).
  \end{split}  
\end{equation*}
To obtain the coefficient of $\tilde t^2$ in $z=\log(y_1(t))$, we consider the
Taylor expansion of $0\equiv G(t,z)=1-tV(e^z)$ around $(\rho,\log\tau)$, where
$\tau=y_1(\rho)$. We set $\tilde z=z-\log\tau$ and write $G$ short for
$G(\rho,\log\tau)$, subscripts denote the partial derivative:
\begin{equation*}
  \begin{split}
    0&=G(t,z)\\
    &=G
    -G_t\tilde t^2
    +G_z\tilde z
    -G_{t,z}\tilde t^2 \tilde z
    +\frac{1}{2}G_{z,z}\tilde z^2
    -\frac{1}{2}G_{t,z,z}\tilde t^2\tilde z^2
    +\frac{1}{6}G_{z,z,z}\tilde z^3 +\dots,
  \end{split}
\end{equation*}
since $G_{t,t}\equiv 0$. We have
\begin{align*}
    G&=0\\
    G_t&=\eval{-V(e^z)}_{(t=\rho,z=\log\tau)}
                                =-\frac{1}{\rho}\\
    G_z&=\eval{-te^zV^\prime(e^z)}_{(t=\rho,z=\log\tau)}
                                =-\rho\tau V^\prime(\tau)=0\\
    G_{t,t}&=0\\
    G_{t,z}&=\eval{-e^z V^\prime(e^z)}_{(t=\rho,z=\log\tau)}
                                     =0\\
    G_{z,z}&=\eval{-t\left(e^{2z}V^{\prime\prime}(e^z)
                                    +e^z V^\prime(e^z)\right)}_{(t=\rho,z=\log\tau)}
                      =-\rho\tau^2 V^{\prime\prime}(\tau)\\
    G_{z,z,z}&=-\rho\tau^3 V^{\prime\prime\prime}(\tau)
                       -\rho\tau^2 V^{\prime\prime}(\tau).
\end{align*}
Therefore, substituting into the Taylor expansion $\tilde z=\alpha\tilde
t+\beta\tilde t^2+\BigO(\tilde t^3)$ we obtain
\begin{equation*}
  \begin{split}
    0&=-G_t\tilde t^2
    +\frac{1}{2}G_{z,z}\tilde z^2
    +\frac{1}{6}G_{z,z,z}\tilde z^3 +\BigO(\tilde t^4)\\
    &=\left(-G_t
       +\frac{1}{2}G_{z,z}\alpha^2\right)\tilde t^2
    +\left(G_{z,z}\alpha\beta
       +\frac{1}{6}G_{z,z,z}\alpha^3\right)\tilde t^3+O(\tilde t^4).
  \end{split}
\end{equation*}
Thus
\begin{align}
  \notag
  \alpha&=\sqrt{\frac{2G_t}{G_{z,z}}}
          =\frac{1}{\rho\tau}\sqrt{\frac{2}{V^{\prime\prime}(\tau)}}\\
  \label{eq:beta}
  \beta&=-\frac{\alpha^2G_{z,z,z}}{6G_{z,z}}
          =-\frac{\tau V^{\prime\prime\prime}(\tau)+3V^{\prime\prime}(\tau)}
                 {3\left(\rho\tau V^{\prime\prime}(\tau)\right)^2}.
\end{align} 

It is easy to check that $\beta$, i.e., the coefficient of $\tilde t^0$ in the
singular expansion of $t\frac{d}{dt}\log(y_1(t))$ is exactly one half of the
residue of $1/y(1-\rho V(y))$ at $\tau$.

In summary, we have shown the following:
\begin{thm}
  Consider bridges with set of steps \Set V.  Let
  $y_1(t),y_2(t),\dots,y_{min_v}(t)$ be the solutions of the functional
  equation $1-t\sum_{v\in\Set V}y^v=0$ that tend to zero as $t$ goes to zero,
  $y_1(t)$ being the branch with real positive Taylor coefficients. Let $\tau$
  be the unique positive real number with $V^\prime(\tau)=0$ and let
  $\rho=\frac{1}{V(\tau)}$. Furthermore, let $\tau_k=y_k(\rho)$ for
  $k\in\{2,3,..,min_v\}$.
  
  Then the constant term in the singular expansion of the generating function
  for bridges is
  \begin{equation}
    \label{eq:constantterm}
    \frac{1}{2}\res_{y=\tau}\frac{1}{y(1-\rho V(y))}
    +\sum_{k=2}^{min_v} \res_{y=\tau_k}\frac{1}{y(1-\rho V(y))}.
  \end{equation}
\end{thm}

Unfortunately, we were not able to show that this expression does not vanish if
the step set contains a step of height strictly greater than one. However, we
have the following conjecture, that we are able to prove partially for some
special cases:
\begin{cnj}
  Let $V(y)$ be a Laurent-polynomial with positive coefficients with highest
  exponent equal to $max_v$ and lowest exponent equal to $-min_v$. Let $\tau$ be
  the unique positive solution of $V^\prime(\tau)=0$ and
  \begin{equation*}
    \begin{split}
      f_<(y)&=\prod_{\substack{V(\kappa)=V(\tau)\\\size{\kappa}<\tau}}(y-\kappa)
             =\sum_{k=0}^{min_v-1} a_k z^k\\
      \intertext{and}
      f_>(y)&=\prod_{\substack{V(\kappa)=V(\tau)\\\size{\kappa}>\tau}}(y-\kappa)
             =\sum_{k=0}^{max_v-1} b_k z^k.
    \end{split}
  \end{equation*}
  Consider the decomposition
  \begin{equation}
    \label{eq:partial}
    \frac{1}{y(V(y)-V(\tau))}=\frac{\alpha+\beta y}{(y-\tau)^2}
    +\frac{p_<(y)}{f_<(y)}+\frac{p_>(y)}{f_>(y)},    
  \end{equation}
  where the degree of $p_<$ is $min_v-2$ and the degree of $p_>$ is $max_v-2$.
  
  Then 
  \begin{align}
    \label{eq:f<}
    0&<a_0<a_1<\dots<a_{min_v-1}\\
    \label{eq:f>}
    0&<b_{max_v-1}<b_{max_v-2}<\dots<b_0,
  \end{align}
  the leading term of $p_<$ is negative if $min_v>1$ and the leading term of
  $p_>$ is positive if $max_v>1$.
\end{cnj}

The negativity of the leading term of $p_<$ in the conjecture would already
imply that the constant term of the singular expansion of $\bar B(t)$ around
$\rho$ does not vanish: $-\frac{1}{\rho}(\frac{\beta}{2}+[y^{min_v-2}]p_<(y))$
is exactly the value given by \eqref{eq:constantterm}. Since replacing $y$ by
$1/y$ in the Laurent-polynomial $V(y)$ changes the sign of $\beta$ in the
decomposition~\eqref{eq:partial}, we can assume that $\beta$ is negative or
zero. In fact, if $min_v=1$, it follows from \eqref{eq:beta} that $\beta$ is
negative. Since $[y^{min_v-2}]p_<(y)<0$ for $min_v>1$, the claim follows.

We can prove parts of the conjecture for $\min_v\leq 3$: In general, it is easy
to see that the product of $f_<$ and $f_>$ has positive coefficients and that
both of their constant terms must be positive.

If $min_v\leq 2$ we can show \eqref{eq:f<} and \eqref{eq:f>} by inductive
arguments. We were also able to check the case $min_v=3$ and $max_v\leq 4$.

If $min_v=2$ we can also show that $[y^{min_v-2}]p_<(y)<0$: in this case,
$p_<(y)$ is constant and equals $1/(\tau_2 V^\prime(\tau_2))$, where $\tau_2$ is
the only negative zero of $V(y)=V(\tau)$ which is smaller than $\tau$ in
modulus. Since $V(y)$ tends to infinity as $y$ approaches $0-$, we have that
$V^\prime(\tau_2)>0$, which implies the claim.

Finally, if $min_v=2$ and the coefficients of $V$ are either zero or one, we can
also conclude that $\beta$ is negative: in this case the numerator of
\eqref{eq:beta} equals 
\begin{equation*}
 -\tau
  V^{\prime\prime\prime}(\tau)-3V^{\prime\prime}(\tau)+3\tau^{-1}V^\prime(\tau)
 =-3 a_{-1}\tau^{-3}+3a_1\tau^{-1}-15a_3\tau+\dots
\end{equation*}
If $a_1=0$ then the above expression is trivially negative. Otherwise we have
to show that $3\tau^{-1}<15\tau$. We show that $\tau>\frac{1}{2}$, which is
sufficient: We have
\begin{align*}
  V^\prime(y)&=-2y^{-3}+\sum_{k\geq -1}k a_k y^{k-1}\\
             &\leq -2y^{-3}+\frac{1}{(1-y)^2}
\end{align*}
which is negative for $y\leq\frac{1}{2}$.

\section{Transcendence}
It is now a simple matter to complete the proof of the main
Theorem~\ref{thm:transcendental}: Since in the circumstances of the theorem the
asymptotic expansion of
$$[t^n x^i] \log B(x;t)=[x^i] \left(H(x)\right)^n [t^n]\log \bar B(t)$$
contains
a term $n^{-2}$, the series $[x^i] \log B(x;t)$ cannot be algebraic.  When $i$
is minimal such that there is at least one walk in the slit plane with steps in
$\mathfrak S$ ending at $(i,0)$, Theorem~\ref{thm:cycle} gives that $[x^i] \log
B(x;t)$ is the generating for such walks. To settle the transcendence of
$S_{i,0}(t)$ for general $i$, we only need to note that $[t^n] \log B(x;t)\sim
c_0\rho^{-n}n^{-1}$, where, as we proved in the last section,
$c_0=\sqrt{\frac{2G_t}{G_{z,z}}}=\frac{1}{\rho\tau_k}\sqrt{\frac{2}{V^{\prime\prime}(\tau_k)}}$,
and thus does not vanish. Hence, the leading term of $[t^nx^i] \log B(x;t)$
contains a factor of $1/n^{3/2}$. Thus, in the convolution formula for
$S_{i,0}(t)$, the term $1/n^2$ in the asymptotic expansion of $[x^i] \log
B(x;t)$ cannot be cancelled by terms of the asymptotic expansion of the product
of two or more functions $S_{i_j,0}$.

The proof of the non-D-finiteness of the other functions can be copied verbatim
from the proof of Proposition~22, page 282 of \cite{BousquetMelou2001}.

\section{Acknowledgements}
I would like to thank Michael Drmota, Bernhard Gittenberger, Bernhard Lamel and
Bodo La\ss\ for numerous stimulating discussions concerning the nature of the
solutions of the functional equation~\eqref{eq:characteristiccurve}. Also I'm
very grateful for two anonymous referees who pointed out numerous mistakes and a
wrong conjecture appearing in the manuscript. And, of course, I would like to
thank Mireille Bousquet-M\'elou for introducing me to the problem and for a
wonderful stay in Bordeaux.
\bibliographystyle{amsplain}

\begin{thebibliography}{99}

\bibitem{BanderierFlajolet2002}
Cyril Banderier and Philippe Flajolet, \emph{Basic analytic combinatorics of
  directed lattice paths}, Theoret. Comput. Sci. \textbf{281} (2002), no.~1-2,
  37--80, Selected papers in honour of Maurice Nivat. \MR{2003g:05006}

\bibitem{BousquetMelou2001}
Mireille Bousquet-M{\'e}lou, \emph{Walks on the slit plane: other approaches},
  Advances in Applied Mathematics \textbf{27} (2001), no.~2-3, 243--288,
  Special issue in honor of Dominique Foata's 65th birthday (Philadelphia, PA,
  2000). \MR{2002j:60076}

\bibitem{Drmota1994}
Michael Drmota, \emph{A bivariate asymptotic expansion of coefficients of
  powers of generating functions}, European Journal of Combinatorics
  \textbf{15} (1994), no.~2, 139--152. \MR{94k:05014}

\bibitem{Flajolet1987}
Philippe Flajolet, \emph{Analytic models and ambiguity of context-free
  languages}, Theoretical Computer Science \textbf{49} (1987), no.~2-3,
  283--309, Twelfth international colloquium on automata, languages and
  programming (Nafplion, 1985). \MR{89e:68067}

\bibitem{FlajoletSedgewick6}
Philippe Flajolet and Robert Sedgewick, \emph{The average case analysis of
  algorithms}, 1994.

\bibitem{Hille}
Einar Hille, \emph{{Analytic function theory. Vol. I, II. 2nd ed. corrected.}},
  {Chelsea Publishing Company}, 1973.

\bibitem{Marden}
Morris Marden, \emph{Geometry of polynomials}, Second edition. Mathematical
  Surveys, No. 3, American Mathematical Society, Providence, R.I., 1966. \MR{37
  \#1562}

\end{thebibliography}

\providecommand{\bysame}{\leavevmode\hbox to3em{\hrulefill}\thinspace}
\providecommand{\MR}{\relax\ifhmode\unskip\space\fi MR }
\providecommand{\MRhref}[2]{%
  \href{http://www.ams.org/mathscinet-getitem?mr=#1}{#2}
}
\providecommand{\href}[2]{#2}

\end{document}